\documentclass[11pt]{amsart} 
\usepackage{amssymb} 
\theoremstyle{plain}
\newtheorem{thm}{Theorem} 
\theoremstyle{definition}
 
\newtheorem{conc}{Conclusion}
\newtheorem{prop}{Proposition}[section]
\newtheorem{cor}[prop]{Corollary} 
\newtheorem{lem}[prop]{Lemma}
\newtheorem{defn}[prop]{Definition} 
\newtheorem{claim}[prop]{Claim}
\theoremstyle{remark} 

\newtheorem*{problem}{Problem}
 
\DeclareMathOperator{\ban}{Ban}
\DeclareMathOperator{\dom}{dom} 
\DeclareMathOperator{\cl}{cl}
\DeclareMathOperator{\cf}{cf}

\DeclareMathOperator{\bad}{Bad}

\DeclareMathOperator{\ran}{ran}

\renewcommand{\[}{\bigl{[}}
\renewcommand{\]}{\bigr{]}}
\newcommand{\sk}{\vskip.05in}
\newcommand{\id}{\mathcal{I}} 
\newcommand{\B}{\mathcal{B}}
\newcommand{\X}{\mathcal{X}}

\newcommand{\az}{{\aleph_0}} 
\newcommand{\ao}{{\aleph_1}}

\newcommand{\restr}{\upharpoonright} 
\newcommand{\forces}{\Vdash}

\newcommand{\uf}{\mathcal{U}}
 
\newcommand{\subs}{\subseteq}
\newcommand{\sups}{\supseteq} 

\renewcommand{\frak}{\mathfrak} 

\newcommand{\omal}{[\omega_1]^{\aleph_0}} 
\newcommand{\I}{\mathbb{I}}
\renewcommand{\a}{\alpha}
\newcommand{\A}{\aleph}
\newcommand{\w}{\omega}
\renewcommand{\b}{\beta}
\numberwithin{equation}{section}
\begin{document}
\title{Gently killing S--spaces} 
\author[T. Eisworth, P. Nyikos,
S. Shelah]{Todd Eisworth, Peter Nyikos, and Saharon Shelah}
\address[Todd Eisworth]{Institute of Mathematics\\ The Hebrew
University\\ Jerusalem 91904\\ Israel}
\email{eisworth@math.huji.ac.il} \address[Peter Nyikos]{Department of
Mathematics\\ University of South Carolina\\ Columbia, SC\\ USA}
\email{nyikos@math.sc.edu} \address[Saharon Shelah]{Institute of
Mathematics\\ The Hebrew University\\ Jerusalem 91904\\ Israel}
\address[Saharon Shelah] 
{Department of Mathematics\\ Rutgers University\\ New Brunswick, NJ
08903\\ USA} \email{shelah@math.huji.ac.il} \keywords{} \subjclass{}
\date{\today}
\thanks{The research of the second author was partially supported by
NSF Grant DMS-9322613.}
\thanks{The research of the third author was partially supported by
NSF grant DMS--9704477 and the Israel Science Foundation founded
by the Israel Academy of Sciences and Humanities. This is publication
number 690 in the list of the third author.} 
\begin{abstract}
We produce a model of ZFC in which there are no locally compact
first countable S--spaces, and in which $2^\az<2^\ao$. A consequence
of this is that in this model there are no locally compact, separable,
hereditarily normal spaces of size $\ao$, answering a question of the
second author \cite{Ny90}.
\end{abstract}
\maketitle

\section{Introduction and Notation}

In Problem 9 of \cite{Ny90}, Nyikos asks if there is a ZFC example of a
separable, hereditarily normal, locally compact space of cardinality
$\aleph_1$. He notes there that for a negative answer, it suffices to
produce a model of set theory in which there are neither Q--sets nor
locally compact, locally countable, hereditarily normal S--spaces.

We provide such a model in this paper. In fact, in our model
$2^{\aleph_0}<2^{\aleph_1}$ (so in particular there are no Q--sets)
and there are no locally compact, first countable S--spaces at all
(hence no locally compact, locally countable, hereditarily normal
S--spaces).

In fact, we obtain something even more general.  Recall that
an S--space is a regular, hereditarily separable space which
is not hereditarily Lindel\"of.  By switching the ``separable''
and ``Lindel\"of'' we get the definition of an L--space.  A 
simultaneous generalization of hereditarily separarable and
hereditarily Lindel\"of spaces is the class of spaces of countable
spread---those spaces in which every discrete subspace is 
countable.  One of the basic facts in this little corner of
set-theoretic topology is that if a regular space of countable spread
is not hereditarily separable, it contains an L--space, and
if it is not hereditarily Lindel\"of it contains an S--space \cite{Ro84}.
In our model, every locally compact 1st countable space of countable spread
is hereditarily Lindelof; consequently, there are no S--spaces
in locally compact 1st countable spaces of countable spread.
This result, reminiscent of one half of a celebrated 1978 result
of Szentmiklossy, will be discussed further at the end of the paper
in connection with a fifty-year-old problem of M.~Kat\v{e}tov.

These concepts and results have elegant translations in terms
of Boolean algebras via Stone duality.  The Stone space $\mathcal S(A)$ of a 
Boolean algebra $A$ is hereditarily Lindel\"of iff every 
ideal of $A$ is countably generated, and first countable
iff every maximal ideal is countably generated.  
Also, $\mathcal S(A)$ is of countable spread iff every minimal set of
generators for an ideal is countable. (An ideal is
said to be minimally generated
if it has a generating set $D$ such that no member of $D$
is in the ideal generated by the remaining members.)  Hence we now
know that $2^{\aleph_0} < 2^{\aleph_1}$ is consistent with the
following statement: if a Boolean algebra $A$ has the property 
 that every minimal set of generators for an ideal is countable,
and every maximal ideal of $A$ is countably generated, then
every ideal of $A$ is countably generated.  On the other hand,
this statement has long
been known to be incompatible with CH.

Note that there are restrictions on such models. In \cite{JKR76} it is
shown that CH implies the existence of a locally compact first
countable S--space, and in Chapter 2 of \cite{To} this is shown to
follow from the weaker axiom $\mathfrak{b}=\aleph_1$.  Thus the fact
that our model satisfies $\mathfrak{b}=\aleph_2$ is no accident of the
proof --- something along these lines is required.

As far as background goes, we will assume a reasonable familiarity
with topological notions such as filters of closed sets and free
sequences. We also use a lot of set theory --- we will assume that the
reader is used to working with proper notions of forcing.  

Our main tool is the use of totally proper notions of forcing that
satisfy the $\aleph_2$--p.i.c.\  (properness isomorphism
condition). We will take a moment to recall the needed definitions.

\begin{defn}
\hfill
\begin{enumerate}
\item Let $P$ be a notion of forcing, and $N$ a countable elementary
submodel of $H(\lambda)$ for some large regular lambda with $P\in
N$. An {\em $(N, P)$--generic sequence} is a decreasing sequence of
conditions $\{p_n:n\in\omega\}\subs N\cap P$ such that for every dense
open $D\subs P$ in $N$, there is an $n$ with $p_n\in D$.
\item A notion of forcing $P$ is said to be totally proper (also known
as NNR proper) if for every $N$ as above and $p\in N\cap P$, there is
an $(N, P)$--generic sequence $\{p_n:n\in\omega\}$ with $p_0=p$ that
has a lower bound.
\end{enumerate}
\end{defn}

The following claim summarizes the properties of totally proper
notions of forcing that we will need.  The proofs are not difficult,
and they are explicitly worked out in \cite{judyandme} and \cite{Ei1}.

\begin{claim}
Let $P$ be a totally proper notion of forcing.
\begin{enumerate}
\item $P$ adds no new reals; in fact, forcing with $P$ adds no new
countable sequences of elements from the ground model.
\item If $G\subs P$ is generic, then $G$ is countably closed. In fact,
every countable subset of $G$ has a lower bound in $G$.
\end{enumerate}
\end{claim} 

The following definition is from Chapter VII of \cite{pif}.

\begin{defn}
\label{picdef}
$P$ satisfies the $\aleph_2$--p.i.c.\  provided the following holds (for
$\lambda$ a large enough regular cardinal):\hfill\\ If\hfill
\begin{enumerate}
\item $i<j<\aleph_2$ \sk
\item $N_i$ and $N_j$ are countable elementary submodels of
$H(\lambda)$ \sk
\item $i\in N_i$, $j\in N_j$ \sk
\item $N_i\cap \aleph_2\subs j$ \sk
\item $N_i\cap i= N_j\cap j$ \sk
\item $h$ is an isomorphism from $N_i$ onto $N_j$ \sk
\item $h(i)=j$ \sk
\item $h$ is the identity map on $N_i\cap N_j$ \sk
\item $P\in N_i\cap N_j$ \sk
\item $p\in N_i\cap P$
\end{enumerate}
then (letting $\dot G$ be the $P$--name for the generic set) there is
a $q\in P$ such that:
\begin{enumerate}
\setcounter{enumi}{10}
\item $q\forces ``(\forall r\in N_i\cap P)[r\in\dot
G\Longleftrightarrow h(r)\in\dot G]\text{''}$ \sk
\item $q\forces``p\in\dot G\text{''}$ \sk
\item $q$ is $(N_i, P)$--generic.
\end{enumerate}
\end{defn}

Notice that if $N_i$ and $N_j$ are as in the above definition, then
$N_i$ and $N_j$ contain the same hereditarily countable sets. This
follows because $h$ is an isomorphism.  In particular,
$N_i\cap\omega_1$ and $N_j\cap\omega_1$ are the same ordinal. We also
note that in both of the previous two definitions, it does not matter
if we require that  the models under consideration contain a fixed
parameter $x\in H(\lambda)$.

The properties of $\aleph_2$--p.i.c.\  forcings that we utilize will be
spelled out when we build our model in the last section of the paper.

\section{Handling P--ideals}

\begin{defn}
A P--ideal in $\omal$ (the set of all countable subsets
of $\omega_1$) is a set $\mathcal{I}\subseteq\omal$ such
that 
\begin{itemize}

\item if $A$ and $B$ are in $\mathcal{I}$, then so is $A\cup B$
\sk
\item if $A\in\mathcal{I}$ and $B\subseteq A$, then $B\in\mathcal{I}$
\sk
\item if $A\in\mathcal{I}$ and $B=^* A$, then $B\in\mathcal{I}$
\sk
\item if $A_n\in\mathcal{I}$ for each $n\in\omega$, then there is an
$A\in\mathcal{I}$ such that $A_n\subseteq^* A$ for each $n$. 
\end{itemize}
\end{defn}

\begin{defn}
Let $\id$ be a P--ideal in $\omal$ generated by a set of size
$\aleph_1$. A generating sequence for $\id$ is a sequence
$\{A_\alpha:\alpha<\omega_1\}$ such that 
\begin{itemize}
\item $A_\alpha\subseteq\alpha$
\sk
\item if $\alpha<\beta$ then $A_\alpha\subseteq^* A_\beta$
\sk
\item if $A\in\id$, then there is an $\alpha$ with $A\subseteq A_\alpha$.
\end{itemize}
Clearly every such $\id$ has a generating sequence.
\end{defn}

Our goal in this section is (assuming CH holds) to define a notion of
forcing (which we call $P(\I)$) with the property that any P--ideal
$\id\subs\omal$ in the ground model admits an uncountable set $A$ in
the extension satisfying $[A]^{\aleph_0}\subs\id$ or
$[A]^{\aleph_0}\cap\id=\emptyset$.   The partial order we use is a
modification of one of the posets from \cite{EiNy1}, itself a
modification of the notion of forcing used in \cite{uristevo}.

Assume CH, and let $\I=\langle \id_\xi:\xi<\kappa\rangle$ be a
sequence of P--ideals in $\omal$. Let
$\{A^\xi_\alpha:\alpha<\omega_1\}$ be a generating sequence for
$\id_\xi$ (such a sequence exists because CH holds).  The notion of
forcing we define depends on our choice of 
generating sequences, but we abuse notation and call the notion of
forcing $P(\I)$.

\begin{defn}
A promise is a function $f$ such that
\begin{itemize}
\item $\dom f$ is an uncountable subset of $\omega_1$
\item $f(\alpha)$ is a finite subset of $\alpha$
\end{itemize}
\end{defn}

\begin{defn}
\label{pidef}
A condition $p\in P(\I)$ is a pair $(a_p, \Phi_p)$ such that
\begin{enumerate}
\item $a_p$ is a function
\sk
\item $\dom a_p$ is a countable subset of $\kappa\times \omega_1$ 
\sk
\item $\ran a_p\subs 2$
\sk
\item for $\xi<\kappa$, $[p]_\xi:=\{\zeta<\omega_1:a_p(\xi,\zeta)=1\}$
is in $\id_\xi$ (so $[p]_\xi=\emptyset$ for all but countably many
$\xi$)
\sk
\item $\Phi_p$ is a countable collection of pairs $(v, f)$, where
$v\subs\kappa$ is finite and $f$ is a promise.
\end{enumerate}
A condition $q$ extends $p$ if
\begin{enumerate}
\setcounter{enumi}{5}
\sk
\item $a_q\sups a_p$, $\Phi_q\sups\Phi_p$
\sk
\item for $(v, f)\in\Phi_p$,
\begin{equation}
Y(v, f,q, p)=\{\alpha\in\dom f:(\forall\xi\in
v)([q]_\xi\setminus[p]_\xi\subs A^\xi_\alpha\setminus f(\alpha))\}
\end{equation}
is uncountable, and
\begin{equation}
(v, f\restr Y(v, f, q, p))\in\Phi_q.
\end{equation}
\end{enumerate}
\end{defn}
 
The intent of $P(\I)$ is to attempt to adjoin for each $\xi<\kappa$
an uncountable set $A_\xi$ with $[A_\xi]^{\aleph_0}$ contained in
$\id_\xi$. A condition gives us an approximation to $A_\xi$ for
countably many $\xi$, as well as some constraints on future growth
of these approximations.  A pair $(v,f)\in\Phi_p$ puts limits on
how our approximation to $A_\xi$ can grow for the finitely many
$\xi\in v$. It may be that the forcing fails to produce an uncountable
$A_\xi$ for some $\xi$, but we show that we can do so in every
situation where we need it.

\begin{defn}
Let $p$ be a condition in $P(\I)$, let $D$ be a dense open subset of
$P(\I)$, and let $v$ be a finite subset of $\kappa$. An ordinal
$\alpha$ is bad for $(v, p, D)$ if there is an
$F_\alpha\in[\alpha]^{<\aleph_0}$ such that there is no $q\leq p$ in
$D$ with $$[q]_\xi\setminus [p]_\xi\subs A^\xi_\alpha\setminus
F_\alpha$$ for all $\xi\in v$.  Let $\bad(v, p, D)$ be the set of
$\alpha<\omega_1$ that are bad for $(v,p, D)$.
\end{defn}
\begin{prop}
\label{badprop1}
$\bad(v, p, D)$ is countable.
\end{prop}
\begin{proof}
Suppose not. Let $f$ be the function with domain $\bad(v,p,D)$ that
sends $\alpha$ to $F_\alpha$, so $f$ is a promise. Let $r$ be the
condition in $P(\I)$ with $a_r=a_p$, and
$\Phi_r=\Phi_p\cup\{(v, f)\}$.  Clearly $r$ extends $p$. Now let
$q\leq r$ be in $D$. By definition, there are uncountably many
$\alpha\in\dom(f)$ such that if$\xi\in v$ then 
$[q]_\xi\setminus[r]_\xi$ is a subset of $A^\xi_\alpha\setminus
f(\alpha)$. This is a contradiction, as any $\alpha$ is bad for $(v,
p, D)$, yet $q\in D$ and $$[q]_\xi\setminus [p]_\xi\subs
A^\xi_\alpha\setminus f(\alpha)$$ for all $\xi\in v$.
\end{proof}

\begin{thm}
\label{thm1}
$P(\I)$ satisfies the $\aleph_2$--p.i.c.\ 
\end{thm}
\begin{proof}
Let $i$, $j$, $N_i$, $N_j$, $h$, and $p$ be as in Definition
\ref{picdef}.  For $r\in N_i\cap P(\I)$, we define
\begin{equation}
\label{2.3}
r\cup h(r):=( a_r\cup h(a_r), \Phi_r\cup h(\Phi_r)).
\end{equation}
\begin{lem}
\label{lem1}
Assume that $r\in N_i\cap P(\I)$.
\begin{enumerate}
\item $r\cup h(r)\in P(\I)$.
\sk
\item $r\cup h(r)$ extends both $r$ and $h(r)$.
\sk
\item If $s\in N_i\cap P(\I)$ and $r\leq s$, then $r\cup h(r)\leq
s\cup h(s)$.
\end{enumerate}
\end{lem}
\begin{proof}
Left to reader.
\end{proof}

Now let $\delta=N_i\cap\omega_1=N_j\cap\omega_1$, and let
$\{D_n:n\in\omega\}$ enumerate the dense open subsets of $P(\I)$ that
are members of $N_i$. Our goal is to build a decreasing sequence of conditions
$\{p_n:n\in\omega\}$ in $N_i\cap P(\I)$ such that $p_0=p$, $p_{n+1}\in
N_i\cap D_n$, and such that the sequence $\{p_n\cup
h(p_n):n\in\omega\}$ has a lower bound $q$. The next lemma shows that
this will be sufficient.

\begin{lem}
\label{lem2}
Let $\{p_n:n\in\omega\}$ be an $(N_i, P(\I))$--generic sequence.
\begin{enumerate}
\item $\{h(p_n):n\in\omega\}$ is an $(N_j, P(\I))$--generic sequence.
\sk
\item If $\{p_n\cup h(p_n):n\in\omega\}$ has a lower bound $q$, then
$q$ satisfies conditions 11 and 13 of Definition \ref{picdef}.
\end{enumerate}
\end{lem}
\begin{proof}
The first clause follows immediately from the fact that $h$ is an
isomorphism mapping $N_i$ onto $N_j$.  For the second clause, note
\begin{equation}
q\forces\text{``}r\in N_i\cap\dot G\text{''}\Longleftrightarrow r\in
N_i\text{ and }\exists n(p_n\leq r).
\end{equation}
This is because for each $r\in N_i\cap P(\I)$, the set of conditions
that extend $r$ or that are incompatible with $r$ is a dense open
subset of $P(\I)$ that is in $N_i$, and hence for some $n$ either
$p_n$ extends $r$ or $p_n$ incompatible with $r$. Similarly, we have
\begin{equation}
q\forces\text{``}r\in N_j\cap\dot G\text{''}\Longleftrightarrow r\in
N_j\text{ and }\exists n (h(p_n)\leq r).
\end{equation}
Now clause 11 of Definition \ref{picdef} follows easily. Clause 13
holds because the $p_n$'s are an $(N_i, P(\I))$--generic sequence.
\end{proof}

Recall that $\delta=N_i\cap\omega_1=N_j\cap\omega_1$, and let
$\{\gamma_n:n\in\omega\}$ enumerate $N_i\cap\kappa$.  We construct by
induction on $n\in\omega$ objects $p_n$, $F_n$, $q_n$ and $u_n$ such
that
\begin{list}{(\roman{enumi})}{\setlength{\leftmargin}{.5in}\setlength{\rightmargin}{.5in}\usecounter{enumi}}
\item $p_0=p$, $F_0=\emptyset$, $u_0=\emptyset$
 \sk
\item $q_n=p_n\cup h(p_n)$ 
\sk
\item $p_{n+1}\in N_i\cap D_n$
 \sk
\item $F_n$ is a finite subset of $\delta$
 \sk
\item $u_n$ is a finite subset of $N_i\cap\kappa$
 \sk
\item $p_{n+1}\leq p_n$
 \sk
\item $F_{n+1}\sups F_n$
 \sk
\item $u_{n+1}\sups u_n$
 \sk
\item $\{\gamma_m:m<n\}\subs u_n$ 
\sk
\item for $\gamma\in u_{n+1}\cup h(u_{n+1})$,
$[q_{n+1}]_\gamma\setminus [q_n]_\gamma\subs A^\gamma_\delta\setminus
F_{n+1}$
 \sk
\item if $(v,f)\in\Phi_{q_k}$ for some $k$,
then there is a stage $n\geq k$ for which
\begin{align}
v&\subs u_{n+1}\cup h(u_{n+1})
\intertext{and}
\label{needed}
\bigl{\{}\alpha\in Y(v,f,q_n, q_k)&: (\forall\xi\in
v)\bigl{(}A^\xi_\delta\setminus F_{n+1}\subs A^\xi_\alpha\setminus
f(\alpha)\bigr{)}\bigr{\}}
\end{align}
is uncountable.
\end{list}
We assume that we have fixed a bookkeeping system so that at each
stage of the induction we are handed a pair $(v, f)$ from some
earlier $\Phi_{q_k}$ for which we must ensure (xi), and such that
every such $(v, f)$ appearing along the way is treated in this manner.

There is nothing to be done at stage $0$, so assume we have carried
out the induction through stage $n$.  At stage $n+1$, we will be
handed $p_n$, $F_n$, $q_n$, and $u_n$, and our bookkeeping hands us
$(v, f)\in \Phi_{q_k}$ for some $k\leq n$.

To start, we choose $u_{n+1}\sups u_n\cup\{\gamma_n\}$ satisfying (v),
but large enough so that $v\subs u_{n+1}\cup h(u_{n+1})$. This means
that (v), (viii), and (ix) hold.

\begin{claim}
If $f$ is a promise, $B\subs\dom f$ uncountable, $v\subs\kappa$
finite, and $\beta<\omega_1$, then there is a finite
$\bar{F}\subs\beta$ such that
$$\bigl{\{}\alpha\in B:(\forall\xi\in v)
\bigl{(}A^\xi_\beta\setminus\bar{F}\subs A^\xi_\alpha\setminus
f(\alpha)\bigr{)}\bigr{\}}$$ is uncountable. 
\end{claim}
\begin{proof}
Straightforward, by induction on $|v|$.
\end{proof}
(Although the preceding claim has a trivial proof, it does not
generalize to the context of the next section and in some sense this
fact is the reason why the next section is so complicated.)
 
Now apply the preceding claim to $v$, $f$, $Y(v,f,q_n, q_k)$,
$u_{n+1}\cup h(u_{n+1})$, and $\delta$ to get a finite
$\bar{F}\subs\delta$ such that
\begin{equation}
\bigl{\{}\alpha\in Y(v,f,q_n, q_k):(\forall\xi\in u_{n+1}\cup
h(u_{n+1}))\bigl{(} A^\xi_\delta\setminus\bar{F}\subs
A^\xi_\alpha\setminus f(\alpha)\bigr{)}\bigr{\}}
\end{equation}
is uncountable. In particular, our choice of $u_{n+1}$ implies
\begin{equation}
\bigl{\{}\alpha\in Y(v,f,q_n, q_k):(\forall\xi\in
v)\bigl{(}A^\xi_\delta\setminus\bar{F}\subs A^\xi_\alpha\setminus
f(\alpha)\bigr{)}\bigr{\}}
\end{equation}
is uncountable. Now let $F_{n+1}=F_n\cup\bar{F}$. Clearly we have
satisfied (iv) and (vii).

Next, we choose $\beta<\omega_1$ such that
\begin{equation}
\label{star3}
N_i\models\beta\notin\bad(p_n, D_n).
\end{equation}
For each $\gamma\in u_{n+1}\cup h(u_{n+1})$ there is a finite
$G_\gamma\subs\beta$ such that $A^\gamma_\beta\setminus G_\gamma\subs
A^\gamma_\delta\setminus F_{n+1}$, so there is a finite $G\subs\beta$
such that
\begin{equation}
\forall\gamma\in u_{n+1}\cup h(u_{n+1})[A^\gamma_\beta\setminus G\subs
A^\gamma_\delta\setminus F_{n+1}].
\end{equation}
Note that both $\beta$ and $G$ are in $N_i\cap N_j$, and hence are
fixed by $h$. By (\ref{star3}), we can find $p_{n+1}\in N_i$ such that
$p_{n+1}\leq p_n$, $p_{n+1}\in D_n$, and
\begin{equation}
N_i\models(\forall\gamma\in u_{n+1})\bigl{(}[p_{n+1}]_\gamma\setminus
[p_n]_\gamma\subs A^\gamma_\beta\setminus G\bigr{)}.
\end{equation}
Applying $h$, we see that
\begin{equation}
N_j\models (\forall\gamma\in
h(u_{n+1}))\bigl{(}[h(p_{n+1})]_\gamma\setminus [h(p_n)]_\gamma\subs
A^\gamma_\beta\setminus G\bigr{)}.
\end{equation}
Thus
\begin{equation}
(\forall\gamma\in u_{n+1}\cup
h(u_{n+1}))\bigl{(}[q_{n+1}]_\gamma\setminus [q_n]_\gamma\subs
A^\gamma_\beta\setminus G\subs A^\gamma_\delta\setminus
F_{n+1}\bigr{)}.
\end{equation}
Our choice of $p_{n+1}$ (and $q_{n+1}$) satisfies (ii), (iii), (vi), and
(x). Since $\bar{F}\subs F_{n+1}$, we have that (xi) is satisfied for
this particular $(v,f)$.

Now we need to verify that the sequence $\{q_n:n\in\omega\}$ has a
lower bound $q$. To start, we define

\begin{equation}
\begin{split}
a_q=&\bigcup_{n\in\omega}a_{q_n}\\
[q]_\xi=&\bigcup_{n\in\omega}[q_n]_\xi
\end{split}
\end{equation}

\begin{claim}
\hfill
\begin{enumerate}
\item $a_q:(N_i\cup N_j)\cap\kappa\rightarrow 2$
\item If $\xi\in N_i\cap\kappa$, then
$[q]_\xi=\cup\{[p_n]_\xi:n\in\omega\}$. If $\xi\in N_j\cap\kappa$,
then $[q]_\xi=\cup\{[h(p_n)]_\xi: n\in\omega\}$.  \sk
\item $[q]_\xi\in\id_\xi$ for $\xi<\kappa$.
\end{enumerate}
\end{claim}
\begin{proof}[Proof of Claim]
Part 1 of the claim follows because the sequence
$\{p_n:n\in\omega\}$ (resp. $\{h(p_n):n\in\omega\}$) meets every dense
set in $P(\I)$ that is a member of $N_i$ (resp. $N_j$). Part 2 follows
as in the proof of Lemma \ref{lem1}.  For the last part, if $\xi\notin
(N_i\cup N_j)\cap\kappa$ there is nothing to check, so assume $\xi\in
(N_i\cup N_j)\cap\kappa$, and fix $n$
such that $\xi\in\{\gamma_n, h(\gamma_n)\}$. Our construction
guarantees that $[q]_\xi\subs [q_n]_\xi\cup A^\xi_\delta$, and this
latter set is in $\id_\xi$.
\end{proof}

\begin{claim}
If $k\in\omega$ and $(v,f)\in\Phi_{q_k}$, then
\begin{equation}
K(v,f,k):=\{\alpha\in\dom f:(\forall\xi\in v)\bigl{(}[q]_\xi\setminus
[q_k]_\xi\subs A^\xi_\alpha\setminus f(\alpha)\bigr{)}\}
\end{equation}
is uncountable.
\end{claim}
\begin{proof}
Let $n\geq k$ be such that our bookkeeping handed us the promise
$(v,f)$ at stage $n+1$ of the construction. The actions we took at
stage $n+1$ ensure that
\begin{equation}
A:=\{\alpha\in Y(v,f,q_n, q_k): (\forall\xi\in
v)(A^\xi_\delta\setminus F_{n+1}\subs A^\xi_\alpha\setminus
f(\alpha)\bigr{)}\}
\end{equation}
is uncountable.  We claim that $A\subs K(v,f,k)$; to see this fix
$\alpha\in A$, and let $\xi\in v$ be arbitrary. We must verify that
$[q]_\xi\setminus [q_k]_\xi$ is a subset of $A^\xi_\alpha\setminus
f(\alpha)$.

\begin{equation*}
\begin{split}
[q]_\xi\setminus [q_k]_\xi =& ([q]_\xi\setminus
[q_n]_\xi)\cup([q_n]_\xi\setminus [q_k]_\xi)\\ \subs& (\bigcup_{m\geq
n}[q_m]_\xi\setminus [q_n]_\xi)\cup A^\xi_\alpha\setminus f(\alpha)\\
\subs& A^\xi_\delta\setminus F_{n+1}\cup A^\xi_\alpha\setminus
f(\alpha)\\ \subs& A^\xi_\alpha\setminus f(\alpha)
\end{split}
\end{equation*}
Notice that in obtaining the second line, we used that $\alpha\in
Y(v,f,q_n, q_k)$, and to obtain the third line we used requirement
(x) of our construction and the fact that $v\subs u_{n+1}\cup
h(u_{n+1})$.
\end{proof}

Now we define
\begin{equation}
\Phi_q=\bigcup_{n\in\omega}\Phi_{q_n}\cup\bigcup_{n\in\omega}\{(v,
f\restr K(v,f,n)):(v,f)\in \Phi_{q_n}\}
\end{equation}
and $q=(u_q, x_q, \Phi_q)$ is a lower bound for the sequence
$\{q_n:n\in\omega\}$ as desired.
\end{proof}

Notice that in our proof, the only relevant properties of $h$ were
that it is an isomorphism from $N_i$ onto $N_j$ that is the identity
on $N_i\cap N_j$ --- the other requirements from Definition
\ref{picdef} were not used. In particular, our proof goes through in
the case that $h$ is actually the identity map (so $N_i=N_j$).  Thus
we obtain the following.
\begin{thm}
\label{piprop}
$P(\I)$ is totally proper.
\end{thm}

We are still not through, however, as we have not yet verified that
$P(\I)$ lives up to its billing.

\begin{defn}
Let $f$ be a promise and $v\subs\kappa$ finite. For $\xi\in v$, we define a set
$\ban_\xi(v,f)$ by 
\begin{equation}
\beta\in\ban_\xi(v,f)\Longleftrightarrow \{\alpha\in\dom f: \beta\in
A^\xi_\alpha\setminus f(\alpha)\}
\end{equation}
is countable. If $\xi\notin v$ then let $\ban_\xi(v,f)=\emptyset$.
\end{defn}

\begin{prop}
If $\xi<\kappa$, and there is no uncountable $A\subs\omega_1$ with
$[A]^{\aleph_0}\cap\id_\xi=\emptyset$, then $\ban_\xi(v,f)$ is
countable.
\end{prop}
\begin{proof}
We can assume that $\xi\in v$ as otherwise
there is nothing to prove. Our assumption on $\id_\xi$ means that
there is an infinite $B\subs A$ with $B\in\id_\xi$.  For each
$\alpha\in\dom f$, there is a finite set $F_\alpha$ for which
$B\setminus F_\alpha\subs A^\xi_\alpha\setminus f(\alpha)$. Thus there is
a single finite $F$ for which
\begin{equation}
\{\alpha\in\dom f : B\setminus F\subs A^\xi_\alpha\setminus
f(\alpha)\}
\end{equation}
is uncountable. Therefore any member of $B\setminus F$ is not in
$\ban_\xi(v,f)$, a contradiction.
\end{proof}

\begin{prop}
If $\xi<\kappa$ and there is no uncountable $A\subs\omega_1$ with
$[A]^{\aleph_0}\cap\id_\xi=\emptyset$, then for each
$\gamma<\omega_1$, the set of conditions $p$ for which
$[p]_\xi\setminus\gamma$ is non--empty is dense in
$P(\I)$.
\end{prop}
\begin{proof}
Let $\xi$ and $\gamma$ be as in the assumption, and let $p\in P(\I)$
be arbitrary. By the previous proposition,
\begin{equation}
\bigcup\{\ban_\xi(v,f):(v,f)\in\Phi_p\}
\end{equation}
is countable (as $\Phi_p$ is countable), hence there is an
$\alpha>\gamma$ not in $\ban_\xi(v,f)$ for any $(v,f)\in\Phi_p$. It is
straightforward to see that there is a $q\leq p$ with 
$\alpha\in[q]_\xi$. 
\end{proof}

\begin{conc}
Assume CH, and let $\I=\langle I_\xi:\xi<\kappa\rangle$ be a list of
P--ideals in $\omal$. Then there is a totally proper notion of forcing
$P(\I)$, satisfying the $\aleph_2$--p.i.c., so that in the generic
extension, for each $\xi<\kappa$ there is an uncountable
$A_\xi\subs\omega_1$ for which either $[A_\xi]^{\aleph_0}\subs\id_\xi$
or $[A_\xi]^{\aleph_0}\cap\id_\xi=\emptyset$.
\end{conc}
\begin{proof}
We have all the ingredients of the proof already. By Theorems
\ref{thm1} and \ref{piprop}, we know $P(\I)$ is totally proper and
satisfies the $\aleph_2$--p.i.c.\  Fix $\xi<\kappa$, assume that $G\subs
P(\I)$ is generic over $V$, and work for a moment in $V[G]$.  

If in $V$ there is an uncountable $A_\xi$ with
$[A_\xi]^{\aleph_0}\cap\id_\xi=\emptyset$, then $A_\xi$ still has this
property in $V[G]$. (Note that since $P(\I)$ is totally proper, no new
countable subsets of $\omega_1$ are added, so $\id_\xi$ is unchanged
by passing to  $V[G]$.) If no such set exists in $V$, then the set 
\begin{equation}
A_\xi:=\bigcup_{p\in G}[p]_\xi
\end{equation}
is uncountable by the previous proposition, and
$[A_\xi]^{\aleph_0}\subs\id_\xi$ by definition of our forcing notion.
\end{proof}

\section{Handling Relevant Spaces}
Our goal in this section is to build, assuming that CH holds, a totally proper
notion of forcing having the
$\aleph_2$--p.i.c.\  that destroys all first countable, countably
compact, non--compact S--spaces in the ground model. In fact, we do a
little better than this --- if $X$ is a first countable, countably
compact, non--compact regular space with no uncountable free sequences, then
after we force with our poset,  $X$ acquires an uncountable free
sequence. The partial order we use is a modification of that used in
\cite{Ei1}, although things do not work as smoothly as they did in the
last section.

\begin{prop}
\label{firstprop}
Suppose $\uf$ is a countably complete (not necessarily maximal) filter
of closed subsets of the 
space $X$, and suppose $Z\subs X$ meets every set in $\uf$. If $\cl_X
Z_0\notin\uf$ for every countable subset $Z_0$ of $Z$, then $X$ has an
uncountable free sequence.
\end{prop}

The proof of the proposition is straightforward. As a corollary, we
note that $\uf$ is generated by separable sets if $X$ has no
uncountable free sequences, and so under CH the filter $\uf$ is
generated by a family of size at most $\aleph_1$.

\begin{defn}
If $\uf$ if a filter of closed subsets of $X$, we say that $Y\subs X$
is $\uf$--large if $Y\cap A\neq\emptyset$ for every $A\in\uf$. We say
that $Y\subs X$
{\em diagonalizes} $\uf$ if $Y$ is $\uf$--large and $Y\setminus A$ is
countable for every set $A\in\uf$.
\end{defn}

Notice that if $\uf$ is countably complete and $\uf$ is generated by a
set of size at most $\aleph_1$, then every $\uf$--large set $Y$ has a
subset $Z$ that diagonalizes $\uf$. If in addition $\uf$ is not fixed,
then every uncountable subset of $Z$ will diagonalize $\uf$ as well.

Let us call a space $X$ {\em relevant} if $X$ is first countable,
countably compact, non--compact, regular, $|X|=\aleph_1$, and $X$ has no
uncountable free sequences. For each relevant $X$, we fix a maximal
filter of closed sets $\uf_X$ that is not fixed. Since we have assumed
CH holds and $X$ is relevant, we can fix a set $Y_X\subs X$ that
diagonalizes $\uf_X$. By passing to a subset if necessary, we may
assume that $Y_X$ is right--separated in type $\omega_1$.

Since $\uf_X$ is a maximal filter of closed sets, this means that
$Y_X$ is a sub--Ostaszewski subspace of $X$, i.e., every closed subset
of $Y_X$ is either countable or co--countable. The filter $\uf_X$ is
reconstructible from $Y_X$ as the set of all closed subsets of $X$
that meet $Y_X$ uncountably often.

We assume that each $Y_X$ has $\omega_1$ as an underlying set, and
that this correspondence is set up so that initial segments are
open. Thus given a collection of relevant spaces, a countable ordinal
$\alpha$ is viewed as a point in each of the spaces.

We also fix a function $\B$ so that for each relevant space $X$ and
ordinal $\alpha<\omega_1$, $\{\B(X,\alpha, n):n\in\omega\}$ is a
decreasing neighborhood base for $\alpha$ as a point in $X$. We will
need one more definition before defining our notion of forcing.

\begin{defn}
A promise $f$ is a function whose domain is an uncountable subset of
$\omega_1$ and whose range is a subset of $\omega$.
\end{defn}

\noindent Until said otherwise, $\X=\{X_\xi:\xi<\kappa\}$ is a
collection of relevant spaces, and CH holds.

\begin{defn}
\label{pxdef}
A condition $p\in P(\X)$  is a pair $(a_p, \Phi_p)$ such that 
\begin{enumerate}
\item $a_p$ is a function
\sk
\item $\dom a_p$ is a countable subset of $\{(\xi,x):\xi<\kappa\text
{ and }  x\in X_\xi\}$
\sk
\item $\ran a_p\subs 2$
\sk
\item for each $\xi<\kappa$, $[p]_\xi:=\{x\in X_\xi: a_p(\xi, x)=1\}$
satisfies $\cl_{X_\xi}[p]_\xi\notin \uf_\xi$
\sk
\item $\Phi_p$ is a countable set of pairs $(v, f)$ where $v\subs
\kappa$ is finite and $f$ is a promise.
\end{enumerate}
A condition $q$ extends $p$ if
\begin{enumerate}
\setcounter{enumi}{5}
\item $a_q\sups a_p$, $\Phi_q\sups \Phi_p$
\sk
\item for $(v, f)\in\Phi_p$, 
\begin{equation}
Y(v,f,q,p):=\{\alpha\in\dom f: (\forall\xi\in v)\[\;[q]_\xi\setminus
[p]_\xi\subs\B(X_\xi,\alpha, f(\alpha))\]\}
\end{equation}
is uncountable, and 
\begin{equation}
(v, f\restr Y(v,f,q,p))\in\Phi_q.
\end{equation}
\end{enumerate}
\end{defn}

The notion of forcing we have described (seemingly) need not be
proper. If, however, we put restrictions on the family $\X$ we get a
proper notion of forcing. We will need some notation to express the
necessary ideas.

\begin{defn}
Let $v\subs\kappa$ be finite. We define
\begin{equation}
X_v=\prod_{\xi\in v}X_\xi,
\end{equation}
and we let $\uf_v$ be the filter of closed subsets of $X_v$ that is
generated by sets of the form $\prod_{\xi\in v}A_\xi$, where $A_\xi\in
\uf_\xi$.   
\end{defn} 

Note that $\uf_v$ will be countably complete and generated by
$\leq\aleph_1$ sets because each $\uf_\xi$ is.

\begin{defn}
Let $v\subs\kappa$ be finite, and let $f$ be a promise. A point
$(x_\xi:\xi\in v)\in X_v$ is banned by $(v, f)$ if
\begin{equation}
\{\alpha\in\dom f: (\forall\xi\in v)\[ x_\xi\in\B(X_\xi,\alpha, f(\alpha))\]\}
\end{equation}
is countable. We let $\ban (v, f)$ be the collection of all points in
$X_v$ that are banned by $(v, f)$.
\end{defn}

\begin{defn}
Let $v\subs\kappa$ be finite. We
say $v$ is {\em dangerous} if there is a promise $f$ such that $\ban
(v, f)$ is $\uf_v$--large. $\X$ is {\em safe} if no finite
$v\subs\kappa$ is dangerous.
\end{defn}

Our definition of ``safe'' was formulated so that the proof of the
following theorem goes through --- the proof of Claim \ref{notban} is
the place where we really need it.

\begin{thm}
\label{thmA}
If $\X=\{X_\xi:\xi<\kappa\}$ is safe, then $P(\X)$ is totally proper.
\end{thm}

Before we commence with the proof of this theorem, we need a
definition and lemma.

\begin{defn}
Let $v\subs\kappa$ be finite, $p\in P(\X)$, and let $D\subs P(\X)$ be
dense. An ordinal $\gamma<\omega_1$ is said to be {\em bad for
}$(v,p,D)$ if there is an $n$ such that there is no $q\leq p$ in $D$
such that for all $\xi\in v$,
\begin{equation}
[q]_\xi\setminus [p]_\xi\subs\B(X_\xi, \gamma, n).
\end{equation}
We let $\bad(v,p,D)$ be the collection of all $\gamma<\omega_1$ that
are bad for $(v,p,D)$.
\end{defn}

So $\gamma\notin\bad(v,p,D)$ means for every $n$, we can find a $q\leq
p$ in $D$ such that $[q]_\xi\setminus [p]_\xi\subs \B(X_\xi,\gamma,
n)$ for all $\xi\in v$.

\begin{lem}
\label{badsmall}
$\bad(v, p, D)$ is countable.
\end{lem}
\begin{proof}
Suppose not. The function $f$ with domain $\bad(v, p, D)$ that sends
$\gamma$ to the $n$ that witnesses $\gamma\in \bad(v,p,D)$ is a
promise. Now we define  $r=(a_p,\Phi_p\cup\{(v, f)\})$. Clearly $r\leq p$ in
$P(\X)$, and since $D$ is dense there is a $q\leq r$ in $D$. Now $Y(v,
f, q, r)$ is uncountable, and for $\gamma\in Y(v,f,q,r)$ and $\xi\in
v$ we have
\begin{equation}
[q]_\xi\setminus [p]_\xi=[q]_\xi\setminus [r]_\xi\subs \B(X_\xi,
\gamma, f(\gamma))
\end{equation}
and this contradicts the definition of $f$.
\end{proof}

\begin{lem}
\label{banclosed}
Let $(v, f)$ be a promise, and suppose $(x_\xi:\xi\in v)$ is not in
$\ban(v, f)$. Then there is $(U_\xi:\xi\in v)$ such that $U_\xi$ is a
neighborhood of $x_\xi\in X_\xi$ and
\begin{equation}
\{\alpha\in\dom f:(\forall \xi\in v)\[ U_\xi\subs\B(X_\xi,\alpha,
f(\alpha))\]\}
\end{equation}
is uncountable. In particular, $\ban (v, f)$ is a closed subset of
$X_v$.
\end{lem}
\begin{proof}
Let $\{V_n:n\in\omega\}$ be a neighborhood base for $(x_\xi:\xi\in
v)$ in the (first countable) space $X_v$, and define
\begin{equation}
A=\{\alpha\in\dom f:(\forall\xi\in v)\[ x_\xi\in\B(X_\xi,\alpha,
f(\alpha))\]\}.
\end{equation}
By assumption, $A$ is uncountable, and for each $\alpha\in A$ there is
an $n$ for which
\begin{equation}
V_n\subs\prod_{\xi\in v}\B(X_\xi,\alpha,f(\alpha)).
\end{equation}
Thus there is a single $n$ for which
\begin{equation}
\{\alpha\in A: V_n\subs\prod_{\xi\in v}\B(X_\xi,\alpha, f(\alpha))\}
\end{equation}
is uncountable. The definition of the product topology then gives us
the $U_\xi$'s that we need.
\end{proof}

\begin{proof}[Proof of Theorem \ref{thmA}]
Let $N\prec H(\lambda)$ be countable with $P(\X)\in N$. Let $p\in
N\cap P(\X)$ be arbitrary, and let $\{D_n:n\in\omega\}$ list the dense
open subsets of $P(\X)$ that are members of $N$. Le
$\delta=N\cap\omega_1$, and let $\{\gamma_n:n<\omega\}$ enumerate
$N\cap\kappa$. 

Since all the spaces in $\X$ are countably compact and $N$ is
countable, there is a sequence $\{\delta_n:n\in\omega\}$ increasing
and cofinal in $\delta$ such that for every $\xi\in N\cap\kappa$, the
sequence $\{\delta_n:n\in\omega\}$ converges in $X_\xi$ to a point
$z_\xi$.

\begin{claim}
\label{notban}
If $v=\{\xi_0,\dots\xi_{n-1}\}\subs N\cap\kappa$ and $f\in N$ is a
promise, then $(z_{\xi_0},\dots,z_{\xi_{n-1}})$ is not banned by $(v,
f)$.
\end{claim}
\begin{proof}
Since $\X$ is safe and $(v, f)\in N$, there are sets
$A_i\in\uf_{\xi_i}\cap N$ for $i<n$ such that $A_0\times\cdots\times
A_{n-1}$ is disjoint to $\ban (v, f)$. Since $A_i\cap\omega_1$ is
countable in $X_{\xi_i}$, for all sufficiently large $\ell$ we have
$\delta_\ell\in A_i$. Since this holds for each $i$, for all
sufficiently large $\ell$ the $n$--tuple
$(\delta_\ell,\dots,\delta_\ell)$ is in $A_0\times\cdots\times
A_{n-1}$. Since this latter set is closed, we have that
$(z_{\xi_0},\dots,z_{\xi_{n-1}})$ is in $A_0\times\cdots\times
A_{n-1}$, hence $(z_{\xi_0},\dots,z_{\xi_{n-1}})$ is not banned by
$(v,f)$.
\end{proof}

Let $\{V(z_\xi, n):n\in\omega\}$ be a decreasing neighborhood base for
$z_\xi$ in $X_\xi$, with $\cl_{X_\xi}V(z_\xi, 0)\notin\uf_\xi$; this
uses the fact that each $X_\xi$ is regular.

We define $p_n\in P(\X)$, $u_n\subs \kappa$ and $g\in\omal$ such that
\begin{enumerate}
\item $p_0=p$, $u_0=\emptyset$, $g(0)=0$
\sk
\item $p_{n+1}\leq p_n$
\sk
\item $p_{n+1}\in N\cap D_n$
\sk
\item $u_n$ is finite
\sk
\item $u_{n+1}\sups u_n$
\sk
\item $h(n+1)>h(n)$
\sk
\item $\{\gamma_m:m<n\}\subs u_n$
\sk
\item for $\gamma\in u_{n+1}$, $[p_{n+1}]_\gamma\setminus
[p_n]_\gamma\subs V(z_\gamma, g(n+1))$
\sk
\item \label{cond} if $(v, f)$ appears in $\Phi_{p_k}$ for some $k$,
then there is an $n\geq k$ for which $v\subs u_{n+1}$ and 
\begin{equation}
\{\alpha\in Y(v,f,p_n,p_k):(\forall\xi\in u_{n+1})\[ V(z_\xi,
g(n+1))\subs\B(X_\xi,\alpha, f(\alpha))\]\}
\end{equation}
is uncountable.
\end{enumerate}

Assume that a suitable bookkeeping procedure has been set up so that
at each stage $n+1$ we are handed a $(v, f)$ in $\Phi_{p_k}$ for some
earlier $k$ for the purposes of ensuring condition \ref{cond}, and in
such a way that every such $(v, f)$ so appears.

There is nothing to be done at stage $0$. At stage $n+1$ we will be
handed $p_n$, $u_n$, and  $g\restr n+1$, and our bookkeeping hands us
a $(v,f)\in\Phi_{p_k}$ for some $k\leq n$.

Choose $u_{n+1}\subs N\cap\kappa$ finite with $u_n\cup v\cup
\{\gamma_n\}\subs u_{n+1}$. Clearly $u_{n+1}$ satisfies 4, 5, and 7.

Let $f'$ be the promise $f\restr Y(v,f,p_n,p_k)$. Clearly $f'$ is in
$N$. By Claim \ref{notban}, we know that
$(z_\xi:\xi\in u_{n+1})$ is not banned by $(u_{n+1}, f')$. Thus by an
application of Lemma \ref{banclosed} we can choose
a value for $g(n+1)>g(n)$ large enough so that
\begin{equation}
\{\alpha\in\dom f':(\forall\xi\in u_{n+1})\[V(z_\xi, g(n+1))\subs
\B(X_\xi,\alpha, f(\alpha))\]\}
\end{equation}
is uncountable.
Now we choose $\ell<\omega$ large enough so that
$\delta_\ell\notin\bad(u_{n+1}, p_n, D_n)$ and
\begin{equation}
(\forall\xi\in u_{n+1})\[\delta_\ell\in V(z_\xi, h(n+1))\]
\end{equation}

Next choose $m$ large enough so that 
\begin{equation}
(\forall\xi\in u_{n+1})\[ \B(X_\xi,\delta_\ell, m)\subs V(z_\xi,
h(n+1))\].
\end{equation}
Since $\B\in N$, we can apply the definition of
$\delta_\ell\notin\bad(u_{n+1},p_n, D_n)$ to get $p_{n+1}\leq p_n$ in
$N\cap D_n$ such that
\begin{equation}
(\forall\xi\in
u_{n+1})\[\;[p_{n+1}]_\xi\setminus[p_n]_\xi\subs\B(X_\xi,\delta_\ell,m)\subs
V(z_\xi, h(n+1))\]. 
\end{equation}
\vskip.2in
\noindent Now why does the sequence $\{p_n:n\in\omega\}$ have a lower bound?

Define $a_q=\bigcup_{n\in\omega}a_{p_n}$ Note that $a_q$ is a function
satisfying requirements 1--3 of Definition \ref{pxdef}, and
$[a_q]_\xi\neq\emptyset$ only if $\xi\in N\cap\kappa$. If $\xi\in
N\cap \kappa$, then $\xi=\gamma_m$ for some $m\in\omega$, and our
construction guarantees that
\begin{equation}
[a_q]_\xi\subs [p_m]_\xi\cup V(z_\xi, 0)
\end{equation}
and so $\cl_{X_\xi}[a_q]_\xi\notin\uf_\xi$.

Now suppose $(v, f)\in\Phi_{p_k}$ for some $k\in\omega$. Define
\begin{equation}
K(v, f, k)=\{\alpha\in\dom f:(\forall\xi\in v)\[\;[x_q]_\xi\setminus
[p_k]_\xi\subs\B(X_\xi,\alpha, f(\alpha))\]\}.
\end{equation}

\begin{claim}
$K(v, f, k)$ is uncountable.
\end{claim}
\begin{proof}
Let $n\geq k$ be as in condition \ref{cond} for $(v, f)$, so
\begin{equation}
A:=\{\alpha\in Y(v, f, p_n, p_k):(\forall\xi\in v)\[ V(z_\xi,
h(n+1))\subs\B(X_\xi, \alpha, f(\alpha))\]\}
\end{equation}
is uncountable.  For $\alpha\in A$ and $\xi\in v$, we have
\begin{equation*}
\begin{split}
[a_q]_\xi\setminus [p_k]_\xi &=\bigcup_{m\geq n}[p_m]_\xi\setminus
[p_n]_\xi\cup [p_n]_\xi\setminus [p_k]_\xi\\
&\subs\bigcup_{m\geq n}[p_m]_\xi\setminus [p_n]_\xi\cup \B(X_\xi,
\alpha, f(\alpha))\quad \text{ (as $A\subs Y(v, f, p_n, p_k)$)}\\
&\subs V(z_\xi, h(n+1))\cup\B(X_\xi,\alpha, f(\alpha))\quad \text{(by
 8 of our construction)}\\
&\subs \B(X_\xi,\alpha, f(\alpha))\qquad\qquad\qquad\text{ (as $\alpha\in A$)}
\end{split}
\end{equation*}
Thus $A\subs K(v, f, k)$.
\end{proof}

So if we define
\begin{equation}
\Phi_q=\bigcup_{n\in\omega}\Phi_{p_n}\cup\bigcup_{n\in\omega}\{(v,
f\restr K(v, f, n)): (v, f)\in\Phi_{p_n}\}
\end{equation}
we have $q=(a_q, \Phi_q)$ is a lower bound for $\{p_n:n\in\omega\}$.
\end{proof}

\begin{prop}
A singleton is safe, so if $\X=\{X\}$ then $P(\X)$ is totally proper.
\end{prop}
\begin{proof}
Suppose $(\{X\}, f)$ form a counterexample. Then $\ban (\{X\}, f)$ is a
$\uf_ X$--large subset of $X$. Since $X$ has no uncountable free
sequences, there is a countable $A=\{x_n:n\in\omega\}\subs\ban(\{X\} ,f)$
such that $\cl_X A\in\uf_X$ and hence
\begin{equation}
B:=\dom f\cap\cl_X A
\end{equation}
is uncountable. If $\alpha\in B$, then there is an $n\in\omega$ with
$x_n\in\B(X,\alpha, f(\alpha))$. Thus there is a single $n$ for which
the set of $\alpha\in B$ with $x_n\in\B(X,\alpha, f(\alpha))$ is
uncountable, and this contradicts the fact that $x_n\in\ban (v, f)$.
\end{proof}

Since the union of an increasing chain of safe collections is itself
safe, we know that maximal safe collections of relevant spaces exist.

\begin{prop}
Assume $\X=\{X_\xi:\xi<\kappa\}$ is safe, $u\subs\kappa$ is finite,
and $p\in P(\X)$. There is a set $A\in\uf_u$ such that for any
$(x_\xi:\xi\in u)\in A$, there is a $q\leq p$ such that
$x_\xi\in[q]_\xi$ for all $\xi\in u$. 
\end{prop}

\begin{proof}
For each $\xi\in u$ we define a set $A_\xi\in\uf_\xi$ as follows:

Let $\{(v_n, f_n):n\in\omega\}$ list all members of $\Phi_p$ with
$\xi\in v_n$ (the assumption that this set is infinite is purely for
notational convenience). For each $n\in\omega$ there is a set
\begin{equation}
B_n:=\prod_{\zeta\in v_n}B^n_\zeta\in\uf_{v_n}
\end{equation}
that is disjoint to $\ban (v_n, f_n)$. Note that this means that for
every $w\subs v_n$ and $(x_\zeta:\zeta\in w)\in\prod_{\zeta\in
w}B^n_\zeta$, the set
\begin{equation}
\label{star2}
\{\alpha\in\dom f_n:(\forall\zeta\in w)\[ x_\zeta\in\B(X_\zeta,\alpha,
f(\alpha))\]\}
\end{equation}
is uncountable.

We let $A_\xi=\bigcup_{n\in\omega}B^n_\xi$, and we check that
$A=\prod_{\xi\in u}A_\xi$ is as required.

So suppose $x_\xi\in A_\xi$ for $\xi\in u$, and define
\begin{equation}
a_q=a_p\cup\{\langle\xi, x_\xi, 1\rangle:\xi\in u\}.
\end{equation}
We want to show that for $(v, f)\in\Phi_p$ the set
\begin{equation}
K(v, f, p)=\{\alpha\in\dom f:(\forall\xi\in v)\[\;
[a_q]_\xi\setminus[p]_\xi\subs\B(X_\xi,\alpha,f(\alpha))\]\}
\end{equation}
is uncountable. Note that this reduces to showing
\begin{equation}
\{\alpha\in\dom f:(\forall \xi\in u\cap v)\[ x_\xi\in\B(X_\xi, \alpha,
f(\alpha))\]\}
\end{equation}
is uncountable, and this follows easily from the fact that the set in
(\ref{star2}) is uncountable.

Thus if we define
\begin{equation}
\Phi_q=\Phi_p\cup\{(v, f\restr K(v, f, p)):(v, f)\in\Phi_p\},
\end{equation}
then $q=(a_q, \Phi_q)$ is the desired extension of $p$.
\end{proof}

\begin{cor}
\label{densesetcor}
If $v\subs\kappa$ is finite, $Z\subs X_v$ is $\uf_v$--large, and $p\in
P(\X)$, then there is a $q\leq p$ and $(x_\xi:\xi\in v)\in Z$ such
that $x_\xi\in [q]_\xi$ for all $\xi\in v$.
\end{cor}

\begin{thm}
\label{sweatthm}
Suppose $\X$ is a maximal safe family, and let $X$ be an arbitrary
relevant space. If $G\subs P(\X)$ is generic, then
\begin{equation}
V[G]\models\text{``$X$ has an uncountable free sequence''.}
\end{equation}
\end{thm}
\begin{proof}
{\noindent CASE 1: $X\in\X$}\hfill\\

In this case $X=X_\xi$ for some $\xi<\kappa$. Let
\begin{equation}
A=\bigcup_{p\in G}[p]_\xi.
\end{equation}
The filter $\uf_\xi$ generates a countably complete filter of closed
subsets of $X_\xi$ in the extension; we will abuse notation a little
bit and call this filter $\uf_\xi$ as well. Note that a set is
$\uf_\xi$--large in $V[G]$ if and only if it meets every set $A\in
\uf_\xi\cap V$.

Now let $A=\bigcup_{p\in G} [p]_\xi$. Clearly $A$ is a subset of
$X_\xi$ in the extension, and since $G$ is countably closed, if we are
given a countable $A_0\subs A$ there is a $p\in G$ with $A_0\subs
[p]_\xi$. This means (in $V[G]$) that the closure of every countable
subset of $A$ is not in $\uf_\xi$.  Given a set $Z\in \uf_\xi$, we can
apply Corollary \ref{densesetcor} with $v=\{\xi\}$ to conclude that
$A\cap Z$ is non--empty. Thus in $V[G]$ the set $A$ is
$\uf_\xi$--large. By Proposition \ref{firstprop}, $X_\xi$ has an
uncountable free sequence.
\vskip.2in
{\noindent CASE 2: $X\notin\X$}\hfill\\
In this case, by the maximality of $\X$ there is a finite
$v\subs\kappa$ such that $\{X_\xi:\xi\in v\}\cup\{X\}$ is
dangerous. To save ourselves from notational headaches, we assume that
$v=n$, and we will refer to $X$ as $X_n$. We will also let $w$ stand
for $n+1$ so the notation $\uf_w$ and $X_w$ will have the obvious
meaning.

Let $f$ be a promise witnessing that $\{X_i:i\leq n\}$ is
dangerous. In $V[G]$, for $i<n$ we let $A_i=\bigcup_{r\in G}[r]_i$ be
the subset of $X_i$ obtained from the generic filter.

By a density argument, there is a $p\in G$ such that $(v,
f)\in\Phi_p$. Thus if $q\leq p$ in $P(\X)$ the set
\begin{equation}
Y(v,f,q,p)=\{\alpha\in\dom f:(\forall
i<n)\[\;[q]_i\setminus[p]_i\subs\B(X_i,\alpha, f(\alpha))\]\}
\end{equation}
is uncountable.

\begin{claim}
\label{important}
In $V[G]$, if $A_i'$ is a countable subset of $A_i\setminus [p]_i$ for
each $i<n$, then 
\begin{equation}
\{\alpha\in\dom f:(\forall i<n)\[
A_i'\subs\B(X_i,\alpha,f(\alpha))\]\}
\end{equation}
is uncountable.
\end{claim}
\begin{proof}
Since $G$ is countably closed, there is a $q\leq p$ in $G$ such that
$A_i'\subs [q]_i\setminus [p]_i$ for all $i<n$. Now we apply the fact
that $Y(v,f,q,p)$ is uncountable.
\end{proof}

Now back in $V$, our assumption is that $\ban(w,f)$ is
$\uf_w$--large. Since $\uf_w$ is $\aleph_1$--complete and generated by
$\aleph_1$ sets, we can choose
\begin{equation}
Z:=\{(x^\xi_i:i\in w):\xi<\omega_1\}\subs\ban(w, f)
\end{equation}
diagonalizing $\uf_w$. By passing to a subsequence, we may assume that
\begin{equation}
\xi_0\neq\xi_1\Rightarrow x^{\xi_0}_i\neq x^{\xi_1}_i
\end{equation}
for all $i\leq n$.  Note also that
\begin{itemize}
\item $\{(x^\xi_i:i<n):\xi<\omega_1\}$ diagonalizes $\uf_v$
\item $\{x^\xi_n:\xi<\omega_1\}$ diagonalizes $\uf_X$
\end{itemize}

\begin{claim}
In $V[G]$, $I=\{\xi<\omega_1:(\forall i<n)x^\xi_i\in A_i\}$ is
uncountable.
\end{claim}
\begin{proof}
This will follow by an easy density argument in $V$. Given $\xi_0<\omega_1$,
the set $\{(x^\xi_i:i<n):\xi\geq\xi_0\}$ still diagonalizes $\uf_v$,
so in particular it is $\uf_v$--large. Now Corollary \ref{densesetcor}
tells us that the set of conditions forcing the existence of a
$\xi>\xi_0$ such that $(\forall i<n)\[x^\xi_i\in[q]_i\]$ is dense in
$P(\X)$, hence $G$ contains such a condition.
\end{proof}

Since $I$ is uncountable, in $V[G]$ the set $\{x^\xi_n:\xi\in I\}$
will diagonalize $\uf_X$.

\begin{claim}
In $V[G]$, if $I_0\subs I$ is countable, then $\cl_X \{x^\xi_n:\xi\in
I_0\}\notin\uf_X$.
\end{claim}
\begin{proof}
Suppose this fails, so there is a countable $I_0\subs I$ witnessing
it. In particular, all but countably many $\alpha<\omega_1$ are in
$\cl_X\{x^\xi_n:\xi\in I_0\}$. For $i<n$, we define
\begin{equation}
A_i'=\{x_i^\xi:\xi\in I_0\},
\end{equation}
and by Claim \ref{important}, the set
\begin{equation}
B=\{\alpha\in\dom f: (\forall i<n)\[ A_i'\subs \B(X_i,\alpha,
f(\alpha))\]\}
\end{equation}
is uncountable. By throwing away a countable subset of $B$, we can
assume that for all $\alpha\in B$, there is a $\xi\in I_0$ such that
$x^\xi_n\in \B(X_n,\alpha, f(\alpha))$.  Thus there is a single
$\xi\in I_0$ for which
\begin{equation}
\{\alpha\in B: x^\xi_n\in\B(X_n,\alpha, f(\alpha))\}
\end{equation}
is uncountable. Now this contradicts the fact that $(x^\xi_i:i\leq n)$
is in $\ban(w, f)$
\end{proof}

We have shown that in $V[G]$, there is a set that diagonalizes $\uf_X$
with the property that the closure of every countable subset is not in
$\uf_X$. Now Proposition \ref{firstprop} tells is that $X$ has an
uncountable free sequence.
\end{proof}

\begin{thm}
If $\X$ is a safe collection of relevant spaces, then $P(\X)$
satisfies the $\aleph_2$--p.i.c.\ 
\end{thm}

\begin{proof}
Let $i$, $j$, $N_i$, $N_j$, $h$, and $p$ be as in Definition
\ref{picdef}. Just as in the previous section, if $r\in N_i\cap
P(\X)$, we define 
\begin{equation}
r\cup h(r):=(a_r\cup h(a_r), \Phi_r\cup h(\Phi_r)).
\end{equation}

\begin{lem} 
Assume that $r\in N_i\cap P(\X)$.
\begin{enumerate}
\item $r\cup h(r)\in P(\X)$
\item $r\cup h(r)$ extends both $r$ and $h(r)$
\item if $s\in N_i\cap P(\X)$ and $r\leq s$, then $r\cup h(r)\leq
s\cup h(s)$
\end{enumerate}
\end{lem}
\begin{proof}
The proof is essentially the same as the one for Lemma \ref{lem1}.
\end{proof}

Just as in the proof of Theorem \ref{thm1}, it suffices to produce an
$(N_i, P(\X))$--generic sequence $\{p_n:n\in\omega\}$ (with $p_0=p$)
such that $\{p_n\cup h(p_n):n\in\omega\}$ has a lower bound.

Let $\{D_n:n\in\omega\}$ list the dense open subsets of $P(\X)$ that
are members of $N_i$. Let $\delta=N_i\cap\aleph_1=N_j\cap\aleph_1$,
and let $\{\gamma_n:n<\omega\}$ enumerate $N_i\cap\kappa$.  Also fix a
sequence $\{\delta_n:n\in\omega\}$ strictly increasing and cofinal in
$\delta$ such that for each $\xi\in (N_i\cup N_j)\cap\kappa$, the
sequence $\{\delta_n:n\in\omega\}$ converges in $X_\xi$ to a point
$z_\xi$.

\begin{claim}
If $v\subs N_i\cap\kappa$ is finite and $f\in N_i$ is a promise, then
$(z_\xi:\xi\in v)$ is not banned by $(v, f)$. The same holds with
$N_i$ replaced by $N_j$.
\end{claim}

For $\xi\in (N_i\cup N_j)\cap\kappa$, let $\{V(z_\xi, n):n\in\omega\}$
be a decreasing neighborhood base for $z_\xi$ in $X_\xi$, with
$\cl_{X_\xi} V(z_\xi, 0)\notin\uf_\xi$.  We will define $p_n$, $q_n$,
$u_n$, and $g\in\vphantom{B}^\omega\omega$ such that
\begin{enumerate}
\item $p_0=p$, $q_0=p_0\cup h(p_0)$, $u_0=\emptyset$, $g(0)=0$
\item $p_{n+1}\leq p_n$
\item $p_{n+1}\in N_i\cap D_n$
\item $q_n=p_n\cup h(p_n)$
\item $u_n\subs N_i\cap \kappa$ is finite
\item $u_{n+1}\sups u_n$
\item $\{\gamma_m:m<n\}\subs u_n$
\item $g(n+1)>g(n)$
\item for $\gamma\in u_{n+1}\cup h(u_{n+1})$,
$[q_{n+1}]_\gamma\setminus [q_n]_\gamma\subs V(z_\gamma, g(n+1))$
\item if $(v, f)\in \Phi_{q_k}$ for some $k$, then there is a stage
$n\geq k$ for which
\begin{equation}
v\subs u_{n+1}\cup h(u_{n+1})
\end{equation}
and
\begin{equation}
\{\alpha\in Y(v, f, q_n, q_k):(\forall\xi\in v)\[V(z_\xi,
g(n+1))\subs\B(X_\xi, \alpha, f(\alpha))\}
\end{equation}
is uncountable.
\end{enumerate}

Fix a bookkeeping procedure as in the proof of Theorem \ref{thm1}. At
stage $n+1$ we will be handed $p_n$, $q_n$, $u_n$, $g\restr n+1$, and
$(v, f)\in\Phi_{q_k}$ for some $k\leq n$.

Choose $u_{n+1}\subs N_i\cap\kappa$ finite with
$u_n\cup\{\gamma_n\}\subs u_n$ and $v\subs u_{n+1}\cup h(u_{n+1})$  To
define $g(n+1)$, we need to split into cases depending on whether $(v,
f)$ comes from $p_k$ or $h(p_k)$.

\noindent Case 1: $(v, f)\in N_i$

Note that $Y(v, f, q_n, q_k)=Y(v, f, p_n, p_k)$, so $f'=f\restr Y(v,
f, p_n, p_k)$ is a promise in $N_i$. We know $(z_\xi:\xi\in v)$ is not
banned by $(v, f')$, hence there is a value $g(n+1)>g(n)$ large enough
such that
\begin{equation*}
\{\alpha\in\dom f':(\forall\xi\in v)\[V(z_\xi,
g(n+1))\subs\B(X_\xi,\alpha, f(\alpha))\]\}
\end{equation*}
is uncountable.

\noindent Case 2: $(v, f)\in N_j\setminus N_i$

This case is analogous --- we use the fact that $Y(v, f, q_n,
q_k)=Y(v, f, h(p_n), h(p_k))$ is in  $N_j$.

In either case, we have ensured that condition (10) of our
construction is satisfied for $(v, f)$.

Now choose $\ell<\omega$ large enough so that
\begin{equation*}
\delta_\ell\notin\bad(u_{n+1}, p_n, D_n)
\end{equation*}
and
\begin{equation*}
(\forall\xi\in u_{n+1}\cup h(u_{n+1}))\[\delta_\ell\in V(z_\xi, g(n+1))\].
\end{equation*}
Choose $m$ large enough so that
\begin{equation*}
(\forall u_{n+1}\cup h(u_{n+1}))\[\B(X_\xi,\delta_\ell, m)\subs
V(z_\xi, g(n+1))\].
\end{equation*}
In $N_i$, apply the definition of $\delta_\ell\notin\bad(u_{n+1}, p_n,
D_n)$ to get $p_{n+1}\leq p_n$ in $N_i\cap D_n$ such that
\begin{equation}
\label{aha1}
(\forall\xi\in
u_{n+1})([p_{n+1}]_\xi\setminus[p_n]_\xi\subs\B(X_\xi,\delta_\ell,
m)).
\end{equation}
Applying the isomorphism $h$ tells us that
\begin{equation}
\label{aha2}
(\forall\xi\in h(u_{n+1}))([h(p_{n+1})]_\xi\setminus
[h(p_n)]_\xi\subs\B(X_\xi,\delta_\ell, m)).
\end{equation}
The choice of $m$, together with (\ref{aha1}) and (\ref{aha2}), tells
us
\begin{equation}
(\forall\xi\in u_{n+1}\cup h(u_{n+1}))([q_{n+1}]_\xi\setminus
[q_n]_\xi\subs V(z_\xi, g(n+1))).
\end{equation}

Thus we have achieved everything required of us at stage $n+1$. The
verification that $\{q_n:n\in\omega\}$ has a lower bound proceeds just
as in the proof of Theorem \ref{thmA}.
\end{proof}
\begin{conc}
Assume CH holds. There is a totally proper notion of forcing $P(\X)$,
satisfying the $\aleph_2$--p.i.c., such that every relevant space in
the ground model acquires an uncountable free sequence in the generic
extension.
\end{conc}

\section{The Iteration}

We now construct a model of ZFC in which $2^{\aleph_0}<2^{\aleph_1}$
and there are no locally compact first countable S--spaces.  Starting
with a ground model $V$ satisfying $2^{\aleph_0}=\aleph_1$ and
$2^{\aleph_1}=\aleph_{17}$, we will do a countable support iteration
of length $\omega_2$.  

More specifically, let $\mathbb{P}=\langle P_\alpha, \dot
Q_\alpha:\alpha<\omega_2\rangle$ be a countable support iteration
defined by
\begin{itemize}
\item $P_0$ is the trivial poset
\item if $\alpha=\beta+1$, then $V^{P_\alpha}\models\dot
Q_\alpha\text{ is Laver forcing}$
\item if $\alpha$ is a limit ordinal, then $V^{P_\alpha}\models \dot
Q_\alpha=\dot P(\I)*\dot P(\X)$, where
$$V^{P_\alpha}\models \I\text{ is the collection of all P--ideals in
}\omal,$$
and
$$V^{P_\alpha*\dot P(\I)}\models\dot\X\text{ is a maximal safe family
of relevant spaces}.$$
\end{itemize}

We don't actually use much about Laver forcing; the relevant facts we
need are that it is proper, assuming CH it satisfies the
$\aleph_2$--p.i.c.\ (Lemma VIII.2.5 of \cite{pif}), and it adds a real
$r\in\vphantom{b}^\omega\omega$ that eventually majorizes every real
in the ground model.

The point of using the partial orders from sections 2 and 3 is that
they can handle all ``candidates'' from a given groundmodel, instead
of just one at a time. This means that in $\omega_2$ stages we can
catch our tail, even though there are $\aleph_{17}$ ``candidates'' to
worry about at each stage of the iteration.

Having defined our iteration, we arrive at the main theorem of this
paper. 

\begin{thm}
\label{MAIN}
In the model $V[G_{\omega_2}]$, there are no locally compact first
countable S--spaces, and $2^{\aleph_0}<2^{\aleph_1}$. More generally,
every locally compact first countable space of countable spread is
hereditarily Lindel\"{o}f.
\end{thm}

The rest of this section will comprise the proof of this theorem. We
start by noting  that for every $\alpha$, 
$$V^{P_\alpha}\models\dot Q_\alpha\text{ has the $\aleph_2$--p.i.c.\  .}$$
This means
\begin{equation}
\label{blah1}
\alpha<\omega_2\Longrightarrow V^{P_\alpha}\models\text{CH}
\end{equation}
(so in particular the definition of $\dot Q_\alpha$ for limit $\alpha$
makes sense) and
\begin{equation}
\label{blah2}
P_{\omega_2}\text{ has the $\aleph_2$--c.c. .}
\end{equation}

The statement (\ref{blah1}) is just Lemma VIII.2.4 of \cite{pif},
while (\ref{blah2}) is Claim VIII.2.9 from the same source.  Note also
that (\ref{blah2}) together with the fact that we are adding many
Laver reals in the iteration implies 
\begin{equation}
\label{blah3}
V^{P_{\omega_2}}\models \mathfrak{b}=2^{\aleph_0}=\aleph_2\text{ and
}2^{\aleph_1}=\aleph_{17}.
\end{equation}
Thus the cardinal arithmetic in $V^{P_{\omega_2}}$ is as advertised,
and we need only verify that every locally compact 1st countable
space of countable spread is hereditarily Lindel\"of
in $V[G_{\omega_2}]$. We first reduce our task by
showing that it suffices to consider only $X$ with a certain form.

\begin{claim}
If $Z$ is a locally compact space of countable spread which is 
not hereditarily Lindel\"of, then there are $X$, $Y$, and
$\{U_\alpha:\alpha<\omega_1\}$ such that
\begin{itemize}
\item $X$ is a locally compact non-Lindel\"of subspace of $Z$ 
\sk
\item $Y\subs X$ is right separated in type $\omega_1$, witnessed by
open sets $\{U_\alpha:\alpha<\omega_1\}$ 
\sk
\item $X=\bigcup_{\alpha<\omega_1}U_\alpha$
\sk
\item $X=\cl Y$
\sk
\item $\ell(X)=\aleph_1$
\end{itemize}
\end{claim}
\begin{proof}

By a basic lemma \cite{Ro84},  $Z$ has a
 right--separated subspace $Y$ of
cardinality $\A_1$, $\{y_\a: \a  < \w_1\}$, and any such subspace
is hereditarily separable
because $Z$ is of countable spread.  For each $y_\a$
pick an open neighborhood $W_\a$ whose closure is compact and
misses all the later $y_\b$. Every locally compact
space is  Tychonoff, so for each $\a$ there is
a cozero-set neighborhood $V_\a$ of $y_\a$ inside $W_\a$.
Let $V = \bigcup \{V_\a : \a \in \w_1\}$.  Then $V$ is locally
compact, and it is not Lindel\"of because each $V_\a$ contains only
countably many $y_\a$.  In fact, $\ell(V) = \A_1$ 
because we carefully took the union of the $V_\a$ instead
of the union of the $W_\a$, and each $V_\a$ is sigma-compact.
Now it is clear that $X = \cl_V Y$ is as desired.
\end{proof}

We work now in the model $V[G_{\omega_2}]$ and assume for purposes of
contradiction that $Z$ is a locally compact first countable
space of countable spread which is not Lindel\"of. Let $X$ and
$Y$ be
as in the previous
claim.  For each $y_\alpha\in Y$, we choose a neighborhood $V_\alpha$
such that $\cl V_\alpha$ is a compact subset of $U_\alpha$. Let
$A_\alpha=V_\alpha\cap Y\in\omal$.

\begin{claim}
$X$ satisfies Property D, i.e., every countable closed discrete subset
of $X$ expands to a discrete collection of open sets.
\end{claim}
\begin{proof}
This follows from the general result that
every 1st countable regular space $X$ satisfying $\ell(X)< \frak b$
satisfies Property D.  The proof of this is only a minor modification
of the proof of \cite[12.2]{vD84} which was for $|X| < \frak b$ because
van Douwen could not find any use for the added generality given by
$\ell(X)< \frak b$.  However, for the sake of self-containment we
give the proof of this result here.
\smallskip
Let $\ell (X) < \frak b$ and let $D = \{x_n: n \in \w\}$ be a countable
closed discrete subspace of $X$.  Using regularity, let 
$\{U_n : n \in \w\}$  be a family of disjoint open sets such that $x_n
\in U_m$ if and only if $x_n = x_m$.  For each $n$ let $\{B_i^n: i \in \w\}$
be a decreasing local base at $x_n$ such that $B_0^n \subset U_n$.   
Let $U = \bigcup \{U_n : n \in \w\}$ 
and for each $y \in Y = X \setminus U$ let $V_y$ be an open neighborhood
of $y$ whose closure misses $D$, and let $f_y: \w \to \w$ be such that
$B_{f_y(n)}^n$ has closure missing $V_y$ for all $n$.  Let 
$\{V_{y_\a} : \a < \kappa \}$ $(\kappa < \frak b)$ cover $Y$ and, using 
the definition of $\frak b,$
let $f : \w \to \w$ be such that $f_{y_\a} <^* f$ for all $\a$.
In other words, there exists $k \in \w$ such that $f_{y_\a}(n) < f(n)$
for all $n \geq k$.  We then have all of $Y$ covered by open sets
each of which meets at most finitely many of the sets $B_{f(n)}^n$,
which is thus a locally finite collection of disjoint open
sets.  Hence it is a discrete open expansion of $D$, as desired.
\end{proof}

Our assumptions on $X$ imply that $|X|\leq\omega_2$ --- every point in
$X$ is the limit of a sequence from $Y$. We will assume that in fact
$|X|=\aleph_2$ (this is the difficult case) and that the underlying
set of $X$ is $\omega_2$, with $Y=\omega_1\subs X$.

Since $X$ is first countable, we have that $w(X)\leq\aleph_2$, so let
$\mathcal{B}=\{W_\xi:\xi<\omega_2\}$ be a base for $X$. For technical
reasons, we assume $W_\xi = U_\xi$ for $\xi<\omega_1$ with repetitions
allowed in the case $w(X)=\aleph_1$.  Let
$\dot\mathcal{B}$ be a $P_{\omega_2}$--name for $\mathcal{B}$, and let
$N$ be an elementary submodel of $H(\lambda)$ satisfying
\begin{itemize}
\item $|N|=\aleph_1$
\sk
\item $X$, $\mathbb{P}$, $\mathcal{B}$, $\dot\mathcal{B}$,
$\{U_\xi:\xi<\omega_1\}$, $\{V_\xi:\xi<\omega_1\}$, and $G_{\omega_2}$
are in $N$
\sk
\item $N\cap\omega_2=\alpha$ for some $\alpha<\omega_2$
\end{itemize}
(The set of such $N$ is closed and unbounded in
$[H(\lambda)]^{\aleph_1}$.)

For an ordinal $\beta<\omega_2$, define
$\mathcal{B}_\beta:=\{W_\xi\cap\beta:\xi<\beta\}$.

\begin{claim}
\label{reflect1}
With $\alpha$ as above,
\begin{enumerate}
\item $\mathcal{B}_\alpha$ is a base for the topology on $\alpha$ as a
subspace of $X$
\item $\mathcal{B}_\alpha\in V[G_\alpha]$
\end{enumerate}
\end{claim}
\begin{proof}\hfill\\
\noindent 1) Suppose $\beta<\alpha$ and $U\subs X$ is a neighborhood
of $\beta$. Since $X$ is first countable and $\beta\in N$, there is a
neighborhood $U'$ of $\beta$ such that $U'\in N$ and $U'\subs U$.  Now
$$N\models(\exists\gamma<\omega_2)[\beta\in W_\gamma\wedge
W_\gamma\subs U'].$$
Thus there is such a $\gamma<\alpha$ and we are done.

\noindent 2) For each pair
$\bar{\beta}=(\beta_0,\beta_1)\in\alpha\time\alpha$, there is a
condition $p_{\bar\beta}\in G_{\omega_2}$ that decides whether or not
$\beta_1\in W_{\beta_0}$, hence there is such a condition in $N$. Now
the support of $p_{\bar\beta}$ is a countable subset of $\omega_2$
that is in $N$, hence there is a $\gamma<\alpha$ with the support of
$p_{\bar\beta}$ a subset of $\gamma$.  This means to decide whether or
not $\beta_1$ is in $W_{\beta_0}$, we need only $\dot\mathcal{B}$ and
$G_{\omega_2}\restr P_\gamma= G_\gamma$.  Thus $\mathcal{B}_\alpha$
can be recovered from $\dot\mathcal{B}$ and the sequence $\langle
G_\gamma:\gamma<\alpha\rangle$, both of which are in $V[G_\alpha]$.
\end{proof}

Now let $\mathfrak{N}=\langle N_\xi:\xi<\omega_2\rangle$ be a
continuous, increasing $\in$--chain of elementary submodels of
$H(\lambda)$ such that
\begin{itemize}
\item each $N_\xi$ is as in the previous discussion
\item $\langle N_\zeta:\zeta<\xi\rangle\in N_{\zeta+1}$
\item $[\omega_2]^{\aleph_0}\subs\bigcup_{\xi<\omega_2}N_\xi$
\end{itemize}

Now we define a function $F:\omega_2\rightarrow\omega_2$ by letting
$F(\xi)$ equal the least $\zeta$ such that
\begin{equation*}
V[G_\xi]\cap[\xi]^{\aleph_0}\subs N_\zeta
\end{equation*}
and
$$N_\xi\cap [\xi]^{\aleph_0}\subs V[G_\zeta].$$
Note that since both $V[G_\xi]\cap[\xi]^{\aleph_0}$ and $N_\xi\cap
[\xi]^{\aleph_0}$ have cardinality at most $\aleph_1$, the function
$F$ is defined for all $\xi<\omega_2$.

\begin{claim}
\label{app1}
Suppose $\alpha<\omega_2$ has cofinality $\aleph_1$ and is closed
under the function $F$.  Then
$N_\alpha\cap[\alpha]^{\aleph_0}=V[G_\alpha]\cap[\alpha]^{\aleph_0}$.
\end{claim}
\begin{proof}
Suppose first that $A\in[\alpha]^{\aleph_0}\cap V[G_\alpha]$. Then
there is a $\beta$ such that $\sup A<\beta<\alpha$ and $A\in V[G_\beta]$.
Now $F(\beta)<\alpha$ and $A\in N_{F(\beta)}\cap
[\beta]^{\aleph_0}\subs N_\alpha\cap [\alpha]^{\aleph_0}$.

Conversely, suppose $A\in [\alpha]^{\aleph_0}\cap N_\alpha$. Since
$\sup A<\alpha$ and $\alpha$ is a limit ordinal, there is a
$\beta>\sup A$ below $\alpha$ with $A\in N_\beta$.  Then $A\in
V[G_{F(\beta)}]\cap[\beta]^{\aleph_0}\subs V[G_\alpha]\cap
[\alpha]^{\aleph_0}$.
\end{proof}

Let $\alpha_0<\omega_2$ be large enough that
$\{A_\xi:\xi<\omega_1\}\in V[G_{\alpha_0}]$ (the $A_\xi$'s were
defined right before Claim 4.2), and let $\alpha<\omega_2$ satisfy
\begin{enumerate}
\item $\alpha>\alpha_0$
\item $\cf(\alpha)=\aleph_1$
\item $N_\alpha\cap\omega_2=\alpha$
\item
$N_\alpha\cap[\alpha]^{\aleph_0}=V[G_\alpha]\cap[\alpha]^{\aleph_0}$.
\end{enumerate}

Such an $\alpha$ can be found by using the preceding claim, as the set
of ordinals satisfying (3) is closed unbounded in $\omega_2$.

\begin{claim}\hfill\\
$V[G_\alpha]\models \id:=\{B\in\omal:|A_\xi\cap B|<\aleph_0\text{ for
all }\xi<\omega_1\}$ is a P--ideal.
\end{claim}
\begin{proof}
Clearly $\id$ is an ideal (and in $V[G_\alpha]$).  Let
$\{B_n:n\in\omega\}\subs\id$ be given; without loss of generality the
$B_n$'s are pairwise disjoint. Since $\cf(\alpha)=\aleph_1$, there is
a $\beta$ in the interval $(\alpha_0, \alpha)$ such that
$\{B_n:n\in\omega\}\in V[G_\beta]$.  For each $\xi<\omega_1$, define a
function $f_\xi\in\vphantom{b}^\omega\omega$ by
\begin{equation}
f_\xi(n)=1+\max(A_\xi\cap B_n).
\end{equation}
Since $\alpha_0<\beta$, each $f_\xi$ is in $V[G_\beta]$. Now in
$V[G_\alpha]$ there is an $r\in\vphantom{b}^\omega\omega$ dominating
$\{f_\xi:\xi<\omega_1\}$ --- $r$ can be taken to be the Laver real
added at stage $\beta+2<\alpha$.  Now let
\begin{equation}
B:=\bigcup_{n\in\omega}B_n\setminus r(n).
\end{equation}
Clearly $B\in\id$ and $B_n\subs^* B$ for all $n\in\omega$.
\end{proof}

Now let $X_\alpha$ be the topological space with
underlying set $\alpha$ and base given by $\mathcal{B}_\alpha$.
Claim \ref{reflect1} tells us that $X_\alpha\in V[G_\alpha]$, and that
in $V[G_{\omega_2}]$, $X_\alpha$ is a subspace of $X$. We will use
this implicitly throughout the remainder of the section.

\begin{claim}\hfill\\
\label{reflect2}
\begin{enumerate}
\item If $A\in V[G_\alpha]\cap [X_\alpha]^{\aleph_0}$ has a limit
point in $X$, then $A$ has a limit point in $X_\alpha$.
\item $V[G_\alpha]\models X_\alpha\text{ has Property D}$
\end{enumerate}
\end{claim}
\begin{proof}\hfill\\
\noindent 1) Suppose $A\in V[G_\alpha]\cap [X_\alpha]^{\aleph_0}$ has
a limit point in $X$.  Our choice of $\alpha$ and Claim \ref{app1}
together imply that $A\in N_\alpha$, and hence there is a limit point
of $A$ in $N_\alpha$.  This gives us the required limit point for $A$ in
$X_\alpha$.
\sk
\noindent 2) Suppose $D=\{x_n:n\in\omega\}$ is a closed discrete subset
of $X_\alpha$ in $V[G_\alpha]$. By the first part of the Claim, $D$ is
a closed discrete subset of $X$, and by Claim \ref{app1} we know that
$D\in N_\alpha$.  Since $X$ satisfies Property D, $D$ expands to a discrete
collection of open sets, without loss of generality members of our
fixed base $\mathcal{B}$. Since $D\in N_\alpha$, there is such an expansion
in $N_\alpha$. Now the countable subset of $\omega_2$ that indexes
this cover is in $N_\alpha\cap[\alpha]^{\aleph_0}$, hence it is in
$V[G_\alpha]$ as well. This gives us the  required discrete family of open sets in $V[G_\alpha]$.
\end{proof}

Our goal is to show that in $V[G_{\alpha+1}]$, $X_\alpha$ acquires an
uncountable discrete subset. Since $X_\alpha$ is a subspace of
$X$ in $V[G_{\omega_2}]$, if we attain our goal we will have a
contradiction, proving that such a space $X$ does not exist in
$V[G_{\omega_2}]$.

We work for a bit in $V[G_\alpha]$.  The first thing we do is force
with $P(\I)$, where $\I$ lists all the P--ideals in $V[G_\alpha]$. If
$H_0$ is a generic subset of $P(\I)$, then in $V[G_\alpha][H_0]$,
either there is an uncountable $B\subs \omega_1$ with
$[B]^{\aleph_0}\subs \id$, or there is an uncountable $B\subs\omega_1$
with $[B]^{\aleph_0}\cap\id = \emptyset$. 

Let us suppose the first possibility occurs. This means that every
countable subset of $B$ has finite intersection with every
$A_\xi$ (in $V[G_\alpha][H_0]$). This continues to hold in
$V[G_{\omega_2}]$, so in $V[G_{\omega_2}]$ there is an uncountable
$B\subs Y$ that meets each $V_\xi$ at most finitely often, i.e., $B$
has no limit points in $Y$. Thus $B$ is a discrete subspace of $Y\subs
X_\alpha$, and we achieve our goal and reach a contradiction.

Now suppose the second possibility occurs. This means that in
$V[G_\alpha][H_0]$, there is an uncountable $B$ such that every
countably infinite subset of $B$ meets some $A_\xi$ in an infinite
set.

\begin{claim}\hfill\\
$V[G_\alpha][H_0]\models Z:=\cl_{X_\alpha} B\text{ is countably
compact and non--compact}$
\end{claim}
\begin{proof}
First note that any countable subset of $Z$ from $V[G_\alpha][H_0]$ is
in $V[G_\alpha]$, as $P(\I)$ is totally proper.  Given $B_0\in
[B]^{\aleph_0}$, there is a $\xi<\omega_1$ such that $B_1=B_0\cap
A_\xi$ is infinite.

Now step into the model $V[G_{\omega_2}]$. Since $B_1\subs A_\xi\subs
V_\xi$ and $\cl V_\xi$ is compact, $B_0$ has a limit point.  Since
$B_0$ is in the model $V[G_\alpha]$, our choice of $\alpha$ implies
that $B_0$ has a limit point in $X_\alpha$.

Now $X_\alpha$ has Property D in $V[G_\alpha]$, and since no new
countable subsets of $X_\alpha$ appear in $V[G_\alpha][H_0]$,
$X_\alpha$ has Property D in this model as well.

This means that any alleged infinite closed discrete subset of $\cl_{X_\alpha}
B$ (in $V[G_\alpha][H_0]$) would expand to a discrete collection of
open sets, thereby yielding an infinite subset of $B$ with no limit
point in $X_\alpha$. We have already argued that this is
impossible. Thus 
$$V[G_\alpha][H_0]\models \cl_{X_\alpha} B\text{ is countably
compact.}$$

Now the open cover $\{X_\alpha\cap U_\xi:\xi<\omega_1\}$ of $X_\alpha$
is in $V[G_\alpha]$ (here we use another assumption we made about
$\mathcal{B}$), and each of these sets meets $B$ at most
countably often, and so $\cl_{X_\alpha} B$ is not compact.
\end{proof}

If it happens that $Z$ contains an uncountable discrete subset, then
we are done, so we may assume this does not happen.  In particular, we
may assume that $Z$ contains no uncountable free sequence. By virtue
of the preceding claim, this means that $Z$ is a relevant space
(terminology from the last section) in $V[G_\alpha][H_0]$.

The next thing we do in our iteration is to force with $P(\X)$, where
$$V[G_\alpha][H_0]\models \X\text{ is a maximal safe collection of
relevant spaces.}$$

The results of the preceding section tell us that $Z$ acquires an
uncountable discrete subset after we do this forcing. Thus

$$V[G_{\alpha+1}]\models X_\alpha\text{ has an uncountable discrete
subset}$$
and again we have achieved our goal, reaching a contradiction. Thus
every first countable locally compact space of countable spread
is hereditarily Lindel\"of; in particular,
there are no locally compact first countable S--spaces in
$V[G_{\omega_2}]$ and Theorem \ref{MAIN} is established.

\bigskip

Theorem 6 is reminiscent of the theorem of Szentmikl\'ossy
recounted in \cite{Ro84} that MA($\w_1$) implies that no compact
space of countable tightness can contain an S--space or an L--space.
Every compact space of countable spread is of countable tightness,
and if a locally compact space is of countable spread, so is
its one-point compactification.  So our result may be looked upon
as a mild version of one half of Szentmikl\'ossy's theorem
for models of $2^{\aleph_0} < 2^{\aleph_1}$.  It would be
very nice if we could get even a similarly mild version of the
other half---it would settle a famous fifty year-old problem of Kat\v{e}tov:

\begin{problem}
 If a compact space has hereditarily
normal (``$T_5$'') square, must it be metrizable?
\end{problem}

The second author showed that the answer is negative if there
is a Q-set, so that in particular MA($\w_1$) implies a negative
answer.  Gary Gruenhage showed that CH also implies a negative
answer.  Proofs appeared in \cite{GN} along with a theorem
connecting Kat\v{e}tov's problem with the theory of S and L spaces:

\begin{thm} 
If there does not exist a Q-set, and  X  is a compact
nonmetrizable space with $T_5$ square, then at least one of the following is
true:
\begin{enumerate} 

\item $X$ is an L-space
\item  $X^2$ is an S-space
\item  $X^2$ is of countable spread, and
 contains both an S-space and an L-space.

\end{enumerate}
\end{thm}

Parts (2) and (3) are ruled out in our model because of Kat\v{e}tov's
theorem that every compact space with $T_5$ square is perfectly
normal, hence first countable.  If it could be shown that there
are no compact L--spaces (which are automatically first countable)
in our model, then Kat\v{e}tov's fifty-year old problem would
be fully solved. It is not out of the question that first countable
compact L--spaces can be gently killed, so that even if some of these
spaces exist in this model, we can maybe throw in a few more notions
of forcing to explicitly banish them.

There is a tantalizing sort of duality between our model and
the model obtained by adding $\aleph_2$ random reals to a model
of MA+$\frak c = \aleph_2$.  There, too, there are no Q-sets
(even though  $2^{\aleph_0} = 2^{\aleph_1}$); but there, it
is L--subspaces of compact spaces of countable spread that
have been ruled out to date, so that (1) and 
(3) that are ruled out there, and it is the status of locally compact first
countable S--spaces that is unknown.

If neither of these models works out, it is to be hoped that
the techniques we have introduced in this paper
will some day produce a model that does settle Kat\v{e}tov's
problem.

\bibliographystyle{amsplain}

\providecommand{\bysame}{\leavevmode\hbox to3em{\hrulefill}\thinspace}

\end{document}